\documentclass[12pt]{article}
\usepackage{amsmath, amssymb,latexsym, epsfig, color}

\newtheorem{theorem}{Theorem}[section]

 \newtheorem{lemma}[theorem]{Lemma}

\def\proof{\smallskip\noindent {\it Proof: \ }}
\def\endproof{\hfill$\square$\medskip}
\def\Z{\mathbb{Z}}

\def\P{\mathcal{ST}}

\DeclareMathOperator{\lk}{lk}
\newcommand{\field}{{\bf k}}

\DeclareMathOperator{\Skel}{Skel}

\DeclareMathOperator{\astar}{ast}

\newcommand{\MON}{\mbox{\upshape \MON}}

\DeclareMathOperator{\Star}{st}

\title{Balanced complexes and complexes without large missing faces}
\author{Michael Goff, Steven Klee, and Isabella Novik
\thanks{Novik's research is partially supported by Alfred P.~Sloan Research
Fellowship and NSF grant DMS-0801152}\\
\small Department of Mathematics, Box 354350\\[-0.8ex]
\small University of Washington, Seattle, WA 98195-4350, USA,\\[-0.8ex]
\small \texttt{[mgoff, klees, novik]@math.washington.edu} }

\begin{document}

\maketitle

\begin{abstract}

The face numbers of simplicial complexes without missing faces
of dimension larger than $i$ are studied.
It is shown that among all such $(d-1)$-dimensional complexes
with non-vanishing top homology, a certain polytopal sphere
 has the componentwise minimal $f$-vector; and moreover, among
all such 2-Cohen--Macaulay (2-CM) complexes, the same sphere
has the componentwise minimal $h$-vector. It is also verified that the
$l$-skeleton of a flag $(d-1)$-dimensional 2-CM complex is
$2(d-l)$-CM while the $l$-skeleton of a flag PL $(d-1)$-sphere
is  $2(d-l)$-homotopy CM.
In addition, tight lower bounds on the face numbers of 2-CM
balanced complexes in terms of their dimension and the number
of vertices are established.
\end{abstract}

\section{Introduction}

In this paper we study balanced simplicial complexes and complexes
without large missing faces. For the latter class of complexes we
settle in the affirmative several open questions raised in the
recent papers by Athanasiadis \cite{Ath} and Nevo \cite{Nevo},
while for the former class we establish tight lower bounds on
their face numbers in terms of dimension and the number of vertices,
thus strengthening the celebrated lower bound theorem for spheres.

A simplicial complex $\Delta$ on the vertex set $[n]:=\{1,2,\ldots, n\}$
is a collection of subsets of $[n]$ that is closed under inclusion and contains
all singletons $\{i\}$ for $i\in[n]$. The elements of $\Delta$ are called its
{\em faces}. A set $F\subseteq [n]$ is called a
{\em missing face} of $\Delta$ if it is not a face of $\Delta$,
but all its proper subsets are. Hence the collection of all missing faces of
$\Delta$ carries the same information as $\Delta$ itself. Thus it
is perhaps not very surprising that imposing certain conditions on
the allowed sizes of missing faces may result in severe restrictions on the
corresponding simplicial complexes.

One simple example of this phenomenon is
that while a simplicial $(d-1)$-sphere may have as few as $d+1$ vertices,
a {\em flag} $(d-1)$-sphere (that is, a simplicial complex with all its
missing faces of size two or, equivalently, 1-dimensional)
needs at least $2d$ vertices. In fact, Meshulam \cite{Mesh} proved that
among all $(d-1)$-dimensional flag simplicial complexes with non-vanishing
top homology, the boundary of the $d$-dimensional cross-polytope simultaneously
minimizes all the face numbers. Similarly, it was recently
verified in \cite{Ath} that among all 2-Cohen--Macaulay (2-CM, for short)
flag (d-1)-dimensional complexes,  the boundary of the $d$-dimensional
cross-polytope simultaneously minimizes all of the $h$-numbers.

In \cite{Nevo}, Nevo considered the more general class of
$(d-1)$-dimensional simplicial complexes with no missing faces of
dimension larger than $i$ (equivalently, of size larger than $i+1$).
He conjectured \cite[Conjecture 1.3]{Nevo} that among all such complexes
 with non-vanishing top homology, a certain polytopal sphere, $S(i, d-1)$
(that for $i=1$ coincides with the boundary of the cross-polytope),
simultaneously minimizes all of the face numbers. He also
asked \cite[Problem 3.1]{Nevo} if the same sphere $S(i, d-1)$
 has the componentwise minimal $h$-vector in the class of all
homology $(d-1)$-spheres without missing faces of dimension larger than $i$.
One of our main results, Theorem~\ref{main-Nevo}, establishes both of these
conjectures.

In addition to verifying that the $h$-numbers of flag spheres
are at least as large as those of the cross-polytope, Athanasiadis shows in
\cite[Theorem 1.1]{Ath} that the graph of a flag simplicial pseudomanifold of
dimension $(d-1)$ is $2(d-1)$-vertex-connected. This is in contrast to the
fact that without the flag assumption one can only guarantee its
$d$-connectedness (for polytopes this is Balinski's theorem, see
\cite[Theorem 3.14]{Zieg95}; the general case is due to Barnette \cite{Bar}).
The above  result prompted Athanasiadis
to ask (see Question~3.2 in the arxiv version of his paper or Remark 3.2
in the journal version) if, for every $0\leq l \leq d-1$,
the $l$-skeleton of a flag homology $(d-1)$-sphere
is $2(d-l)$-CM and if the $l$-skeleton of a flag PL $(d-1)$-sphere
is $2(d-l)$-homotopy CM.  In Theorem~\ref{main-Ath} we settle
both of these questions in the affirmative.

The face numbers of flag complexes are closely related to those of balanced
complexes. (A simplicial $(d-1)$-dimensional complex is called {\em balanced}
\cite{St79} if its 1-skeleton, considered as a graph, is vertex
$d$-colorable.) Indeed, it is a result of Frohmader \cite{Fro}
that for every flag complex $\Delta$ there exists a balanced complex $\Gamma$
with the same $f$-vector, and it is a conjecture of Kalai \cite[p.~100]{St96}
that if $\Delta$ is flag and CM, then one can choose the corresponding
balanced $\Gamma$ to also be CM.

The lower bound theorem for
spheres \cite{Bar-lower, Kalai} asserts that among all homology $(d-1)$-spheres
on $n$ vertices, a stacked sphere has the componentwise minimal $f$-vector.
Here we provide a sharpening of these bounds for the class of
balanced homology spheres in Theorem \ref{balanced}.
In the case of balanced $(d-1)$-spheres
whose number of vertices, $n$, is divisible by $d$, our result amounts
to the statement that the spheres obtained by taking the connected sum
 of $\frac{n}{d} -1$ copies of the boundary
of the $d$-dimensional cross-polytope have the componentwise minimal
$f$-vector.

The rest of the paper is structured as follows. In Section 2 we review basic
facts and definitions related to simplicial complexes and their face numbers.
Section 3 is devoted to complexes without large missing faces. Section 4
deals with CM connectivity of skeletons of flag complexes.
Finally, in Section 5 we discuss balanced complexes.
Sections 3--5 are independent of each other and can be read in any order.
We hope that our results will be helpful
in attacking additional stronger conjectures proposed in \cite{Nevo}.

\section{Preliminaries}

Here we review basic facts and definitions related to simplicial
complexes. An excellent reference to this material is Stanley's
book \cite{St96}.

Let $\Delta$ be a simplicial complex on the vertex set $[n]$.
For $F\in \Delta$, set $\dim F:= |F|-1$ and
define the dimension of $\Delta$, $\dim \Delta$, as the maximal dimension
of its faces. We say that $\Delta$ is \textit{pure} if all of its facets
(maximal faces under inclusion) have the same dimension.
The $f$-vector of $\Delta$ is
$f(\Delta)=(f_{-1}, f_0, \ldots, f_{d-1})$, where $d-1=\dim\Delta$ and $f_j$ is
the number of $j$-dimensional faces of $\Delta$. Thus $f_{-1}=1$ (unless
$\Delta$ is the empty complex)
and $f_0=n$. We also consider the $f$-polynomial of $\Delta$,
$$f(\Delta, x):= \sum_{j=0}^{d} f_{j-1}x^j.$$

It is sometimes more convenient to work with the $h$-vector,
$h(\Delta)=(h_0,h_1,\ldots, h_d)$ (or the $h$-polynomial,
$h(\Delta, x):= \sum_{j=0}^{d} h_{j}x^j$)
instead of the $f$-vector ($f$-polynomial, resp.).
It carries the same information as the $f$-vector and
is defined by the following relation:
$$
h(\Delta, x)= (1-x)^d f\left(\Delta, \frac{x}{1-x}\right).
$$
In particular, $h_0=1$, $h_1=n-d$, and the $f$-numbers of $\Delta$
are non-negative linear combinations of its $h$-numbers.

Let $\Delta_1$ and $\Delta_2$ be simplicial complexes
on disjoint vertex sets $V_1$ and $V_2$. Then their \textit{join}
is the following simplicial complex
on $V_1\cup V_2$,
$$
\Delta_1\ast\Delta_2:=\{F_1\cup F_2 \ : \ F_1\in \Delta_1,  F_2\in \Delta_2\}.
$$
Therefore, $f(\Delta_1\ast\Delta_2, x)=f(\Delta_1, x)f(\Delta_2, x)$
and  $h(\Delta_1\ast\Delta_2, x)=h(\Delta_1, x)h(\Delta_2, x)$.
Also, a set $F\subseteq V_1\cup V_2$
is a missing face of $\Delta_1\ast\Delta_2$ if and only if it is a
missing face of either $\Delta_1$ or $\Delta_2$.
Thus if both complexes have no
missing faces of dimension larger than $i$, then so does their join.

Similarly, if $\Delta_1$ and $\Delta_2$ are pure simplicial $(d-1)$-dimensional
complexes on disjoint vertex sets, and $F_1=\{v_1,\ldots, v_d\}\in\Delta_1$ and
$F_2=\{w_1,\ldots,w_d\}\in\Delta_2$ are facets, then the complex
obtained from $\Delta_1$ and $\Delta_2$ by identifying $F_1$ and $F_2$ via
the bijection $\rho(v_i)=w_i$, and then removing this identified face, is called
the {\em connected sum} of $\Delta_1$ and $\Delta_2$ along $F_1$ and $F_2$,
and is denoted $\Delta_1\#_{\rho} \Delta_2$. While the combinatorics of the
resulting complex depends on $F_1$, $F_2$, and $\rho$, its
$f$- and $h$-vectors do not:
$$
h_i(\Delta_1\# \Delta_2)=\left\{
\begin{array}{ll}
h_i(\Delta_1)+h_i(\Delta_2)-1 & \mbox{ if $i=0$ or $d$}\\
h_i(\Delta_1)+h_i(\Delta_2) & \mbox{ if $0<i<d$.}
\end{array}
\right.
$$

If $\Delta$ is a simplicial complex and $F$ is a face of $\Delta$,
 then the {\em link}
of $F$ in $\Delta$ is
 $\lk_{\Delta} F=\lk F:= \{G\in\Delta \ : \
 F\cup G\in\Delta, \,\, F\cap G =\emptyset \},
$
the {\em star} of $F$ in $\Delta$ is
$\Star_{\Delta} F=\Star F:=\{G\in\Delta \ : \  F\cup G\in\Delta\}$,
and the {\em antistar} of $F$ in $\Delta$ is
$\astar_{\Delta} F=\astar F=\{G\in\Delta \ : \ F\not\subseteq G\}$.
Also, for $W\subseteq [n]$, let
$\Delta_{-W} :=\{F\in\Delta \ : \ F\subseteq [n]-W\}$ denote the
 {\em restriction} of
$\Delta$ to $[n]-W$. The links, stars, antistars, and restrictions
are simplicial complexes in their own right. If $\Delta$ is a complex
without missing faces
of dimension larger than $i$, then so are links, stars, and restrictions
 of $\Delta$;  furthermore this property is preserved
under taking antistars of faces of dimension at most $i$.

We say that a $(d-1)$-dimensional complex $\Delta$ is
{\em Cohen--Macaulay} over $\field$ (CM, for short)
if $\tilde{H}_i(\lk F; \field)=0$ for all $F\in\Delta$ and all $i<d-|F|-1$.
Here $\field$ is either a field or $\Z$ and $\tilde{H}_i(-, \field)$
denotes the $i$th reduced simplicial homology with coefficients in $\field$.
If in addition, $\tilde{H}_{d-|F|-1}(\lk F; \field)\cong \field$ for every
 $F\in \Delta$, then $\Delta$ is a {\em \field-homology sphere}.
We say that $\Delta$ is {\em $q$-CM}
if for all $W\subset [n]$, $|W|\leq q-1$, the complex
$\Delta_{-W}$ is CM  and has the same dimension as $\Delta$. 2-CM complexes
are also known as doubly CM complexes. Every simplicial
 sphere (that is, a simplicial complex whose geometric realization
is homeomorphic to a sphere) is a homology sphere (over any \field),
and every \field-homology sphere is doubly CM over \field.
Moreover, joins and connected sums of (homology) spheres are
(homology) spheres.

Similarly, we say that  $\Delta$ is
{\em homotopy Cohen--Macaulay} (homotopy CM, for short) if $\lk F$ is
$(d-|F|-2)$-connected for all $F\in\Delta$, and that
$\Delta$ is {\em $q$-homotopy CM} if
$\Delta_{-W}$ is homotopy CM  and has the same dimension as $\Delta$
for all $W\subset [n]$, $|W|\leq q-1$.
(Recall that a complex, or more precisely, its geometric realization,
 is {\em $i$-connected} if all of its homotopy groups
from 0th to the $i$th one vanish.)
Unlike the usual Cohen--Macaulayness, homotopy Cohen--Macaulayness is not
a topological property: there exist simplicial spheres that are not homotopy
CM. It is however worth pointing out that all PL simplicial spheres are
homotopy CM (in fact, 2-homotopy CM).

Two simplicial complexes are said to be {\em PL homeomorphic} if there exists
a piecewise linear map between their geometric realizations that is
also a homeomorphism. A simplicial complex
is a {\em PL $(d-1)$-sphere}
if it is PL homeomorphic to the boundary of the $d$-simplex.
The importance of PL spheres is that all their links are also
PL spheres (see e.g.~\cite[Section 12(2)]{Bj}).

\section{Counting face numbers}
The goal of this section is to prove the following result conjectured in
\cite{Nevo}. Throughout this section we fix positive integers
$i$ and $d$ and write $d=qi+r$
where $q$ and $r$ are (uniquely defined) integers satisfying $1\leq r \leq i$.
Let $\sigma^j$ denote the $j$-dimensional simplex, $\partial \sigma^j$ its
boundary complex, and $(\partial \sigma^j)^ {\ast q}$ the join of $q$
copies of $\partial \sigma^j$. Define
\begin{equation}  \label{S(i,d-1)}
 S(i, d-1) :=
 (\partial \sigma^i)^ {\ast q} \ast \partial \sigma^r.
\end{equation}
Note that $S(1,d-1)$ coincides with the boundary of the $d$-dimensional
cross-polytope.

\begin{theorem} \label{main-Nevo}
Let $\Delta$ be a $(d-1)$-dimensional simplicial complex without missing
faces of dimension larger than $i$, and let $\field$ be a field or $\Z$.
\begin{enumerate}
\item
If $\Delta$ has a non-vanishing top homology (with coefficients in \field),
then $f_j(\Delta)\geq f_j(S(i,d-1))$ for all $j$.  If
$f_0(\Delta)=f_0(S(i, d-1))$
and, unless $d$ is divisible by $i$,
$f_r(\Delta)=f_r(S(i, d-1))$, then $\Delta=S(i, d-1)$.

\item If $\Delta$ is 2-CM over \field, then $h_j(\Delta)\geq h_j(S(i,d-1))$
for all $j$.  If $h_1(\Delta)=h_1(S(i, d-1))$
and, unless $d$ is divisible by $i$,
$h_{r+1}(\Delta)=h_{r+1}(S(i, d-1))$, then $\Delta=S(i, d-1)$.
\end{enumerate}
\end{theorem}

\noindent

Several cases of Theorem \ref{main-Nevo} are known:
Nevo \cite[Theorem 1.1]{Nevo} verified the inequalities in Part 1
for all $j$ assuming that $i$ divides $d$, and for all $j\leq r$ if
$i$ does not divide $d$; the $i=1$ case of Part~2 is due to
Athanasiadis \cite[Theorem 1.3]{Ath}.

Throughout the proof, the inequality $P(x)\geq Q(x)$ between two polynomials
means that the polynomial $P(x)-Q(x)$ has non-negative coefficients.
The proof of both parts relies on the following simple property of the
$h$-numbers of $S(i, d-1)$.
\begin{lemma}  \label{ineq}
For every $1\leq s \leq i$, one has
\begin{eqnarray}
 h(S(i, d-1), x) &\leq& h(\partial \sigma^s, x) h(S(i, d-1-s), x), \quad
\mbox{and hence also}  \label{h-ineq}\\
 f(S(i, d-1), x) &\leq& f(\partial \sigma^s, x) f(S(i, d-1-s), x).
\nonumber
\end{eqnarray}
%Moreover, equality (in both of these instances)
%holds if and only if $s=i$ or $s=r$.
\end{lemma}

\proof Since the $f$-numbers are non-negative combinations of the $h$-numbers,
it is enough to verify the first inequality.
Express $d$ as $d=s+(q'i+r')$, where $q', r'$ are integers
satisfying $1\leq r' \leq i$. Then $q-q'\in\{0,1\}$ and $s+r'=(q-q')i+r$. Since
$h(\partial \sigma^j, x)=\sum_{l=0}^j x^l$, the inequality in (\ref{h-ineq})
divided by $(\sum_{l=0}^i x^l)^{q'}$ reads
$$ \left( \sum_{l=0}^i x^l \right)^{q-q'} \left(\sum_{l=0}^r x^l \right) \leq
\left(\sum_{l=0}^{r'} x^l \right)\left( \sum_{l=0}^s x^l \right).$$
If $q=q'$, then $r=r'+s$, and the above inequality holds without equality.
Otherwise, $q-q'=1$ and $i+r=r'+s$ with $i=\max\{r,r',s,i\}$, and the
assertion follows by comparing coefficients.
\endproof

{\medskip\noindent {\it Proof of Theorem \ref{main-Nevo}, Part 1: \ }}
We first prove the inequalities on the $f$-numbers of $\Delta$
by induction on $d$.  If $d\leq i$, then $S(i,d-1)=\partial \sigma^d$,
and the result follows from the well-known and easy-to-prove fact that
among all simplicial complexes of dimension $(d-1)$ with non-vanishing
top homology, $\partial\sigma^d$ has the componentwise minimal $f$-vector.
 So assume that $d>i$ and that the statement holds for all $d'<d$.

If $F$ is a face with $0\leq \dim F \leq i$ and
$\tilde{H}_{d-1-|F|}(\lk F; \field) = 0$, then consider
$\Delta':=\astar F$ and $\Delta'' := \Star F$, so that
$\Delta = \Delta' \cup \Delta''$ and $\Delta' \cap \Delta''
 = \partial \overline{F} \ast \lk F$. (Here $\overline{F}$
denotes the simplex $F$ together with all its faces.)
Since $\tilde{H}_{d-2}(\Delta' \cap \Delta''; \field)
\cong \tilde{H}_{d-1-|F|}(\lk F; \field)=0$ and since $\Delta''$
 is a cone, and hence acyclic, the Mayer-Vietoris sequence \cite[p.~229]{Br}
 yields that $\tilde{H}_{d-1}(\Delta'; \field)
\cong \tilde{H}_{d-1}(\Delta; \field)\neq 0$.
Therefore, by considering $\Delta'$ instead of $\Delta$, we may assume
without loss of generality that every face $F$ of $\Delta$ with
$\dim F \leq i$ satisfies $\tilde{H}_{d-1-|F|}(\lk F; \field) \neq 0$.

Let $G$ be a missing face of $\Delta$ and consider $G' \subsetneq G$. Define
$\Delta^{G'}$ to be the collection of faces of $\Delta$ of the form
$G' \cup F$, where $F \cap G = \emptyset$.  Note that $\Delta^{G'}$
is not generally a simplicial complex.  Since $H \cap G = G'$ for all
$H \in \Delta^{G'}$, the collections $\Delta^{G'}$ are pairwise  disjoint as
$G'$ ranges over all proper subsets of $G$.  For $G'\subset G$, choose
$G' \subseteq G'' \subset G$  satisfying $|G''|=|G|-1$. Since $G$ is a
missing face in $\Delta$, $\lk G''$ does not contain any vertices from  $G$,
 and therefore $F \cup G' \in \Delta^{G'}$ for all $F \in \lk G''$.
Hence $f(\Delta^{G'},x) \geq x^{|G'|}f(\lk G'',x) \geq x^{|G'|}f(S(i,d-|G|),x)$
by the inductive hypothesis. As the collections $\Delta^{G'}$ are pairwise
disjoint for $G' \subsetneq G$, by summing over all such $G'$, we obtain

\begin{eqnarray*}
f(\Delta, x) & \geq &
  \sum_{G' \subsetneq G} x^{|G'|}f(S(i,d-|G|),x)
 =  f(S(i,d-|G|),x)\sum_{G' \subsetneq G} x^{|G'|}    \\
%& = & f(S(i,d-|G|),x)\sum_{j=0}^{|G|-1} f_{j-1}(\partial G)x^{j} \nonumber\\
& = & f(S(i,d-|G|),x)f(\partial \overline{G},x)
 \geq  f(S(i,d-1),x),
\end{eqnarray*}
where the last step is by Lemma \ref{ineq}.

We now prove the statement on equality by induction on $d$.
Assume that $f_0(\Delta) = f_0(S(i,d-1))$ and, if $r<i$, that
$f_r(\Delta) = f_r(S(i,d-1))$.  Then
$f_j(\Delta) = \binom{f_0(\Delta)}{j+1} = f_j(S(i,d-1))$ for all $j < r$.
Furthermore, $\Delta$ has a missing face of dimension $r$.
 In the case that $r < i$ this follows from
$f_r(\Delta)=f_r(S(i, d-1))=\binom{f_0(\Delta)}{r+1}-1$.
When $r=i$, this follows from the fact that $\Delta$ has a complete
$(r-1)$-dimensional skeleton and no missing face of $\Delta$
has dimension greater than $i$.
Finally, $\tilde{H}_{d-1-|G|}(\lk G; \field) \neq 0$ for all
$G\in \Delta$ with $\dim G < r$:
 otherwise $\tilde{H}_{d-1}(\astar G;\field) \neq 0$
and $f_{\dim G}(\astar G) < f_{\dim G}(S(i,d-1))$, a contradiction.

Let $F$ be a missing face of $\Delta$ of dimension $r$ and $G$ a maximal
proper subset of $F$.  We claim that if $F'$ is a missing face in $\lk G$
of dimension $i$, then $F'$ is a missing face in $\Delta$ as well.  Let $G'$
be a minimal subface of $G$ such that $\lk G'$ does not contain $F'$ as a face.
 Then every proper subface of $G' \cup F'$ is a face in $\Delta$, but not
$G' \cup F'$ itself.  Since $\dim F' = i$, we infer that $G' = \emptyset$ and
$F'$ is a missing face in $\Delta$.

We have that $f_0(\lk G) \leq f_0(\Delta)-r-1$, since $\lk G$ contains no
vertex of $F$; and, in fact, equality holds here by the inductive hypothesis
since $\lk G$ has nonvanishing top homology. Also
$\dim(\lk G)+1 = \dim(\Delta)+1-r=d-r$ is divisible by $i$,
and so it follows by the inductive hypothesis that $\lk G = S(i,d-1-r)$.
 Label the missing faces of $\lk G$ by $F_1, \ldots, F_q$.
Every missing face of $\lk G$ has dimension $i$, and hence
 every missing face of
 $\lk G$ is also a missing face of $\Delta$ by the previous paragraph.
Thus $\Delta$ has $F,F_1, \ldots, F_q$ as disjoint missing faces with
$\dim F = r$ and $\dim F_1 = \ldots = \dim F_q = i$.  These are precisely the
missing faces of $S(i,d-1)$, and so $\Delta$ is contained in $S(i,d-1)$. Since
$S(i,d-1)$ has componentwise minimal face numbers, $\Delta = S(i,d-1)$.
\endproof

The proof of Part 2 utilizes
the following results in addition to Lemma \ref{ineq}.
The first of them is due to
Stanley \cite[Cor.~II.3.2]{St96}, the second appears in works of
Adin, Kalai, and Stanley, see e.g. \cite{St93}, and the third one is
\cite[Lemma~4.1]{Ath}.

\begin{lemma} \label{Stanley}
If $\Delta$ is a $(d-1)$-dimensional CM complex, then $h(\Delta, x)\geq 0$.
Moreover, if $\Delta$ has a non-vanishing top homology (which happens, for
instance, if $\Delta$ is a homology sphere, or more generally, if $\Delta$
is 2-CM), then $h(\Delta, x)\geq \sum_{l=0}^d x^l= h(\partial \sigma^d, x)$.
\end{lemma}

\begin{lemma}  \label{Adin-Kalai-Stanley}
Let $\Delta$ be a simplicial complex and $\Gamma$ a subcomplex of $\Delta$.
 If $\Delta$ and $\Gamma$ are both CM (over the same \field)
and have the same dimension, then $h(\Delta, x)\geq h(\Gamma, x)$.
\end{lemma}

\begin{lemma} \label{Christos}
Let $\Delta$ be a pure simplicial complex and $v$ a vertex of $\Delta$.
If $\astar_\Delta v$ has the same dimension as $\Delta$, then
$h(\Delta, x) = x h(\lk_\Delta v, x) + h(\astar_\Delta v, x)$.
\end{lemma}

We are now in a position to prove Part 2 of the theorem. It follows the same
general outline as the proof of \cite[Theorem 1.3]{Ath}, but requires a bit
more bookkeeping.

{\medskip\noindent {\it Proof of Theorem \ref{main-Nevo}, Part 2: \ }}
The proof is by induction on $d$. If $d\leq i$, then
$S(i,d-1)=\partial \sigma^d$, and the statement follows
from Lemma \ref{Stanley}. So assume that $d>i$ and that the statement holds
for all $d'<d$.

Let $F=\{v_0, v_1, \ldots, v_s\}$ be a missing face of $\Delta$ (in particular,
$s\leq i$). Then
$F_j:=\{v_0, v_1, \ldots, v_j\}$ is a face for every $-1\leq j \leq s-1$,
and so is $F-v_j:=\{v_0, \ldots, \hat{v_j}, \ldots, v_s\}$ for every
$0\leq j\leq s$. Repeatedly applying Lemma \ref{Christos} and using the
fact that $\lk_{\lk G}H = \lk_{\Delta}(H \cup G)$ for all $G \in \Delta$ and
$H \in \lk_{\Delta}G$ (here and below, $\lk$ without a
subscript refers to the link in $\Delta$), we obtain
\begin{eqnarray}
h(\Delta, x) & = &
  x h(\lk_\Delta v_0, x) + h(\astar_\Delta v_0, x) \nonumber \\
& = & x(x h(\lk_{\lk v_0} v_1, x) + h(\astar_{\lk v_0} v_1, x))
    + h(\astar_\Delta v_0, x) \nonumber   \\
& = & x^2(x h(\lk_{\lk F_1} v_2, x) + h(\astar_{\lk F_1} v_2, x))
   +x h(\astar_{\lk F_0} v_1, x) +  h(\astar_\Delta v_0, x) \nonumber\\
& = & \cdots \nonumber \\
& = & x^s h(\lk_\Delta F_{s-1}, x) +
\sum_{j=0}^{s-1} x^{j} h(\astar_{\lk F_{j-1}} v_{j}, x).  \label{sum}
\end{eqnarray}

Since $\Delta$ is 2-CM, all its links are also 2-CM \cite{Baclaw},
and so all the complexes
appearing in (\ref{sum}) are CM. We now show that the $h$-polynomial of
each of these complexes is (componentwise) at least as large as $h(S(i, d-s-1),
 x)$, and hence
\begin{eqnarray*}
h(\Delta, x) &\geq& \sum_{j=0}^s x^j  h(S(i, d-s-1), x)=
\left(\sum_{j=0}^s x^j\right) h(S(i, d-s-1), x) \\
& = & h(\partial \sigma^s, x) h(S(i, d-s-1) \geq h(S(i, d-1), x)
\quad \mbox{ (by Lemma \ref{ineq})},
\end{eqnarray*}
as required.

And indeed, $\lk_\Delta F_{s-1}$ is $(d-s-1)$-dimensional, 2-CM, and has
no missing faces of size larger than $i$. Hence
$h(\lk_\Delta F_{s-1}, x)\geq h(S(i, d-s-1), x)$ by the inductive hypothesis.
For all other complexes appearing in (\ref{sum}), observe that since $F$ is a
missing face, the complex
$v_s \ast v_{s-1} \ast \cdots \ast v_{j+1} \ast \lk_\Delta(F-v_j)$
is well-defined, does not contain $v_j$, and is contained in
$\lk_\Delta F_{j-1}$. In other words,
$$\astar_{\lk F_{j-1}} v_j \supseteq
 v_s \ast v_{s-1} \ast \cdots \ast v_{j+1} \ast \lk_\Delta(F-v_j).$$
As both of these complexes are CM of dimension  $d-j-1$,
Lemma \ref{Adin-Kalai-Stanley} yields that
$$
h(\astar_{\lk F_{j-1}} v_j, x)
\geq h( v_s \ast \cdots \ast v_{j+1} \ast \lk (F-v_j), x )=
h( \lk (F-v_j), x)\geq h(S(i, d-s-1), x),$$
where the last step is by the inductive hypothesis. This completes the proof.
The treatment of equality follows from the first part and the observation
that $S(i,d-1)$ has a complete $(r-1)$-dimensional skeleton.
\endproof

\section{Cohen--Macaulay connectivity of flag complexes}
This section is devoted to the proof of the following theorem. Recall that the
$l$-skeleton of a simplicial complex $\Delta$, $\Skel_l(\Delta)$, consists of
all faces of $\Delta$ of dimension at most $l$.

\begin{theorem} \label{main-Ath}
Let $\Delta$ be a flag simplicial complex of dimension $d-1$.
\begin{enumerate}
\item If $\Delta$ is 2-CM over \field, then $\Skel_l(\Delta)$ is
$2(d-l)$-CM over $\field$ for all $0\leq l\leq d-1$.
\item Moreover, if $\Delta$ is a simplicial PL sphere,
then $\Skel_l(\Delta)$ is
$2(d-l)$-homotopy CM for all $0\leq l\leq d-1$.
\end{enumerate}
\end{theorem}

Throughout the proof,
 $\|\Delta\|$ stands for the geometric realization of $\Delta$;
for $W\subset [n]$, $\overline{W}$ denotes the simplex on the vertex set $W$
together with all its faces,
and $p_W$ denotes the barycenter of $\|\overline{W}\|$.
If $\Gamma$ is a subcomplex of $\Delta$, and  $W$ is a subset
of $[n]$ (but not
necessarily a subset of $V(\Gamma)$ --- the vertex set of $\Gamma$),
we write $\Gamma_{-W}$ to denote the restriction of $\Gamma$ to $V(\Gamma)-W$.
We  make use of the following observation: for $F\in \Delta$ and
$W\subseteq [n]-F$,
\begin{equation}   \label{links}
\lk_{\left(\Skel_l(\Delta)\right)_{-W}} F =
\left(\lk_{\Skel_l(\Delta)} F\right)_{-W}
=\left(\Skel_{l-|F|}(\lk_\Delta F)\right)_{-W}.
\end{equation}

{\medskip\noindent {\it Proof of Part 1: \ }}
In the following $\field$ is fixed and is suppressed from our notation.
The proof is by induction on $d$.
Since $\Delta$ is flag and 2-CM, we already
know that it has at least $2d$ vertices, and hence
that $\Skel_0(\Delta)$ is $2d$-CM. This implies the assertion for $d\leq 2$
as well as for $l$=0 and any $d$.

Assume now that the statement holds for all $d'<d$. In particular,
it holds for all links of non-empty faces of $\Delta$
since they are also 2-CM
and have dimension strictly smaller than $d-1$.
Thus for a nonempty face $F\in\Delta$, the complex $\Skel_{l-|F|}(\lk F)$
is $2((d-|F|)-(l-|F|))=2(d-l)$-Cohen--Macaulay.
Putting this together with (\ref{links})
and using that for $j<l$ the $j$th simplicial homology of $\Skel_l(\Delta)$
coincides with that of $\Delta$, to
complete the proof it only remains to show that (i) for every $W\subset [n]$
of size $2(d-l)-1$, $\Delta_{-W}$ is at least
$l$-dimensional, and (ii) for all $j<l\leq d-1$
and any subset
$W=\{v_1, \ldots, v_k\}\subset [n]$ of size $1\leq k\leq 2(d-l)-1$,
the homology $\tilde{H}_j(\Delta_{-W})$ vanishes.

To verify (i) consider $F\in\Delta_{-W}$ of dimension at most
$(l-1)$. We need to show that $F$ is not a maximal (under inclusion)
face in $\Delta_{-W}$. Since the link of $F$ in $\Delta$
is a flag 2-CM complex of dimension $\geq (d-l-1)$,
it has at least $2(d-l)>|W|$ vertices. Thus, at least one of these
vertices, say, $v$ is not in $W$, yielding that
 $F\cup\{v\}\in \Delta_{-W}$ is a larger face.

To prove (ii) we induct on $k$. There are two possible cases to consider.

\smallskip\noindent{\bf Case 1:} every two vertices of $W$ are connected by
an edge in $\Delta$ (this, for instance, happens if $k$=1).
Since $\Delta$ is flag, this condition implies that $W\in\Delta$.
Then $\|\Delta_{-W}\|$ is a strong deformation retract of
$\|\Delta\|-\|\overline{W}\|$ (see e.g. \cite[Lemma 11.15]{Bj})
which in turn is a strong deformation retract of
$\|\Delta\|-p_W$. Since $\Delta$ is 2-CM, the latter complex is
$(d-2)$-acyclic (this is essentially due to Walker, see
\cite[Prop.~III.3.7]{St96}),  and the statement follows.

\smallskip\noindent{\bf Case 2:} not every two vertices of $W$
form an edge. By reordering the vertices, if necessary, assume that
$\{v_{k-1}, v_k\}\notin \Delta$. Consider complexes
$\Delta_{k-1}:=\Delta_{-(W-v_k)}$, $\Star_{\Delta_{k-1}} v_k$,
and the intersection
 $\Delta_{-W}\cap \Star_{\Delta_{k-1}} v_k=\lk_{\Delta_{k-1}} v_k$.
The first two complexes have vanishing $j$th homology: indeed, the star is
contractible and for $\Delta_{k-1}$ this holds by our inductive hypothesis
on $k$. Also, since $v_{k-1}$ and $v_k$ are not connected by an edge,
$v_{k-1}$ is not in the link of $v_k$. Hence
$$
\lk_{\Delta_{k-1}} v_k =
\left(\lk_\Delta v_k\right)_{-\{v_1,\ldots, v_{k-2}\}}.
$$
But $\dim (\lk_{\Delta} v_k)= d-2$ and
$k-2\leq 2(d-l)-3=2((d-1)-l)-1$, so our inductive hypothesis on $d$ applies to
$\lk_{\Delta} v_k$ and shows that
$\tilde{H}_{j-1}(\lk_{\Delta_{k-1}} v_k)=0$. Finally, since
$\Delta_{k-1}=\Delta_{-W}\cup \Star_{\Delta_{k-1}} v_k$, the appropriate
portion of the Mayer--Vietoris sequence \cite[p.~229]{Br} yields that
$\tilde{H}_j(\Delta_{-W})=0$, and the assertion follows. \endproof

We now turn to Part 2 of the theorem. A PL sphere is 2-CM over $\Z$,
so Part 1 implies
vanishing of relevant homology groups computed with coefficients in $\Z$.
In particular, all the spaces involved are (path) connected, and this allows
us to
suppress the base point when discussing homotopy groups. We also write
$\pi_j(\Delta)$ instead of $\pi_j(\|\Delta\|)$.

The Hurewicz theorem \cite[p.~479]{Br}
asserts that if $\Delta$ is $j$-connected, $j\geq 1$,
then $\pi_{j+1}(\Delta) \cong \tilde{H}_{j+1}(\Delta;\Z)$.  In particular,
if $\Delta$ is simply connected and $\tilde{H}_i(\Delta;\mathbb{Z})=0$
for all $0 \leq i \leq j$, then $\tilde{H}_{j+1}(\Delta;\mathbb{Z})
\cong \pi_{j+1}(\Delta).$
Also, PL spheres are simply connected and their links
are PL spheres in their own right.
Thus Part 2 will follow from Part 1 if we can show that for a
PL $(d-1)$-sphere $\Delta$
and an arbitrary $W\subset [n]$ of size $1\leq k\leq 2(d-2)-1=2d-5$,
$\pi_1(\Delta_{-W})=0$. This is done exactly as in the proof of Part 1:
except that in Case 2 one needs to use
the Seifert--van Kampen theorem \cite[p.~161]{Br}
instead of the Mayer--Vietoris sequence.
It asserts (using notation of Case 2 in the proof of Part 1) that
$$
\pi_1(\Delta_{k-1}) \cong \pi_1(\Delta_{-W})
   \ast_{\pi_1(\lk_{\Delta_{k-1}} v_k)} \pi_1(\Star_{\Delta_{k-1}} v_k).
$$
Since by the inductive hypothesis all groups, except possibly
$\pi_1(\Delta_{-W})$, in this equation are trivial,
it follows that $\pi_1(\Delta_{-W})$ is trivial as well. As for Case 1,
just notice that a topological sphere with a point removed is a topological
ball, and hence contractible. \endproof

We close this section with several remarks.

\smallskip\noindent{\bf 1.} In Part 2 of the theorem the
`PL sphere' condition cannot be relaxed to the `triangulated sphere' one.
This can be seen by considering the double suspension of the Poincar\'e sphere.
According to Edwards, see \cite{Dav},
the resulting space is a topological sphere.
Now start with any triangulation of the Poincar\'e sphere, and let $\Gamma$
be its barycentric subdivision.
Then $\Delta=(\partial \sigma^1)^{\ast 2} \ast \Gamma$
is a flag complex that triangulates Edwards' sphere. But $\Delta$
is not homotopy CM: indeed, some of the edges of $\Delta$ have $\Gamma$ as
their link, and $\Gamma$ is not simply connected.

\smallskip\noindent{\bf 2.} Let $n\geq 2d$ be any integer, and
let $C_k$ denote the graph-theoretical cycle on $k$ vertices. Then
the complex $S(1, d-1, n):=(\partial \sigma^1)^{\ast (d-2)}\ast C_{n-2d+4}$
is a $(d-1)$-dimensional polytopal sphere on $n$ vertices. It is flag, and
for all $1\leq l\leq d-1$, its $l$-skeleton is $2(d-l)$-CM, but
not $(2(d-l)+1)$-CM. Thus Part 1 of Theorem \ref{main-Ath} is as strong as
one can hope for. This example together with the theorem
also adds plausibility to Conjecture 1.4
from \cite{Nevo}
asserting that among all flag homology $(d-1)$-spheres on $n$ vertices,
$S(1, d-1, n)$ has the smallest face numbers.

\smallskip\noindent{\bf 3.} An immediate consequence of Part 1 of the theorem
is that if $\Delta$ is a (d-1)-dimensional flag complex that is
$k$-CM for some $k\geq 2$,
then $\Skel_l(\Delta)$ is $(2(d-l-1)+k)$-CM for $0\leq l \leq d-1$.
Indeed to show that $\Skel_l(\Delta)_{-W}$ is CM and of dimension $l$ for
any $|W| \leq 2(d-l)+k-3$, consider a subset $W'$ of $W$ of size
$\min\{k-2, |W|\}$ and its complement $W''=W-W'$. Then $|W''|\leq 2(d-l)-1$.
Since $\Delta_{-W}=(\Delta_{-W'})_{-W''}$, and since $\Delta_{-W'}$ is 2-CM,
Theorem~\ref{main-Ath} applied to $\Delta_{-W'}$ and $W''$ completes the proof.

It is interesting to compare this result with
 a theorem of Fl{\o}ystad \cite{Fl} asserting
that if $\Delta$ is an arbitrary $(d-1)$-dimensional $k$-CM simplicial
complex, then its $l$-skeleton is $((d-l-1)+k)$-CM.

\section{The lower bound theorem for balanced complexes}
In this section we establish tight lower bounds
on the face numbers of balanced 2-CM complexes
in terms of their dimension and the number of vertices.
Recall that a $(d-1)$-dimensional complex $\Delta$
on the vertex set $V$ is {\em (completely) balanced}
if its 1-dimensional skeleton is $d$-colorable:
that is, there exists a coloring
$\kappa: V \rightarrow [d]$ such that for all
$F \in \Delta$ and distinct $v,w \in F$,
$\kappa(v) \neq \kappa(w)$.  We assume that a balanced
complex $\Delta$ comes equipped with such a coloring $\kappa$.
The order complex of a  rank $d$ graded poset
is one example of a balanced simplicial complex.

If $\Delta$ is a balanced complex and $T \subseteq [d]$, then
the {\em $T$-rank selected subcomplex} of $\Delta$ is
$\Delta_T := \{F \in \Delta: \kappa(F) \subseteq T\}$.
We make use of the following basic facts from \cite{St79}.

\begin{lemma}  \label{St-balanced}
Let $\Delta$ be a $(d-1)$-dimensional balanced CM complex.
Then for any $T\subseteq [d]$,
$\Delta_T$ is also CM, and
$h_i(\Delta) = \sum_{|T|=i}h_i(\Delta_T)$ for all $0\leq i \leq d$.
\end{lemma}

\noindent Since deleting a vertex commutes with taking a rank
selected subcomplex: $(\Delta_T)_{-v}=(\Delta_{-v})_T$ for
any $v$ with $\kappa(v)\in T$, one consequence of the
above lemma is that a rank selected subcomplex of a 2-CM complex
is 2-CM as well.

The Lower Bound Theorem for simplicial spheres \cite{Bar-lower, Kalai}
asserts that among all $(d-1)$-dimensional homology
spheres with $n$ vertices, a stacked sphere,
$\P(n, d-1)$, has the componentwise minimal $f$-vector.
A {\em stacked sphere}, $\P(n, d-1)$,
is defined as the connected sum of $n-d$ copies of the boundary of the
$d$-simplex. Since
$h_1(\partial\sigma^d)=h_2(\partial\sigma^d)=1$ if $d\geq 2$, it follows that
for $d\geq 3$, $h_1(\P(n, d-1))=h_2(\P(n, d-1))$.
Therefore, via a well-known reduction due to McMullen, Perles, and Walkup
(see \cite[Thm.~1]{Bar-lower} or \cite[Sect.~5]{Kalai}),
the proof of the LBT for $d\geq 3$ reduces to showing that the $h$-vector of
a homology sphere of dimension at least $2$ satisfies $h_2\geq h_1$.
Recently, Nevo \cite{Nevo-2CM} extended this result to
all 2-CM simplicial complexes:

\begin{lemma} \label{Nevo-lower}
If $\Delta$ is a simplicial 2-CM complex of dimension at least $2$,
then $h_2(\Delta)\geq h_1(\Delta)$.
\end{lemma}

It follows easily from the results of \cite{St79}
that the boundary of the $d$-dimensional cross-polytope
has the componentwise minimal $h$-vector among all balanced
$(d-1)$-dimensional spheres. This motivates us to define
a {\em stacked cross-polytopal sphere}, $\P^{\times}(n, d-1)$, for $n$
a multiple of $d$, as the connected sum of $\frac{n}{d}-1$ copies of
the boundary complex of the $d$-dimensional cross polytope.
At each step the vertices of the same colors are identified to guarantee that
the resulting complex is balanced as well.

What are the $h$-numbers of $\P^{\times}(n, d-1)$? Since the $h$-numbers
of the $d$-dimensional cross-polytope are given by $h_j=\binom{d}{j}$,
it follows that for  $0<j<d$,
$$
h_j(\P^{\times}(n, d-1))=\left(\frac{n}{d}-1\right)\cdot\binom{d}{j},
\mbox{ that is, }  (j+1)h_{j+1}=(d-j)h_j \mbox{ for } 0<j<d-1.
$$
In particular, $(d-1)h_1=2h_2$ if $d \geq 3$.
Similarly, a direct computation shows that
$$
\psi_{j-1}(n,d-1):= j\cdot f_{j-1}(\P^{\times}(n, d-1))=\begin{cases}
(2^j-1)\binom{d-1}{j-1}(n-d) +
d\binom{d-1}{j-1},&  1 \leq j \leq d-1 \\
(2^d-2)(n-d) + 2d, &  j=d.
\end{cases}
$$
One advantage of the last expression is that it is defined for {\bf all} $n$
rather than just multiples of $d$. This allows us to state and prove
the main theorem of this section --- the Lower Bound Theorem
for balanced spheres and, more generally, balanced 2-CM complexes.

\begin{theorem} \label{balanced}
Let $\Delta$ be a balanced 2-CM simplicial complex of dimension $d-1$.
If $d\geq 3$,
then $ 2h_2(\Delta)\geq (d-1)h_1(\Delta)$. In particular, if $d\geq 2$ and
$f_0(\Delta)=n$,
then $j\cdot f_{j-1}(\Delta)\geq \psi_{j-1}(n,d-1)$ for all $2\leq j\leq d$.
\end{theorem}

\proof
Repeatedly applying Lemma \ref{St-balanced}, we see that
\begin{eqnarray*}
\sum_{|T|=3} h_2(\Delta_T)& = &
\sum_{|T|=3}
\sum_{\mbox{\tiny{$\begin{array}{cc}S\subset T,\\
|S|=2\end{array}$}}} h_2(\Delta_S)
= \sum_{|S|=2} (d-2)h_2(\Delta_S) =
(d-2)h_2(\Delta), \mbox{ and } \\
\sum_{|T|=3} h_1(\Delta_T)& = & \sum_{|T|=3}
\sum_{\mbox{\tiny{$\begin{array}{cc}S\subset T,\\
|S|=1\end{array}$}}} h_1(\Delta_S)
= \sum_{|S|=1} \binom{d-1}{2}h_1(\Delta_S) =
\binom{d-1}{2} h_1(\Delta).
\end{eqnarray*}
Since $\Delta$ is balanced and 2-CM, its rank selected subcomplexes share
the same properties. In particular, when $|T|=3$,
$\Delta_T$ is a 2-dimensional 2-CM
complex, and so by Lemma \ref{Nevo-lower}, $h_2(\Delta_T) \geq  h_1(\Delta_T)$.
Thus we infer that $(d-2)h_2(\Delta)\geq \binom{d-1}{2} h_1(\Delta)$, and
the inequality $2h_2(\Delta)\geq (d-1)h_1(\Delta)$ (for $d\geq 3$) follows.

The proof of the ``in particular'' part is a routine computation
similar in spirit to the McMullen-Perles-Walkup reduction. We sketch it
here for completeness. We use induction on $d$. For $d=2$ we need
only show that $2f_1(\Delta)\geq 2n$. This indeed holds, since $\Delta$ is a
2-CM graph, hence it is 2-connected, and so every vertex of $\Delta$
has degree at least 2.

Suppose now that $d\geq 3$. Then
$\sum_{v \in \Delta}h_1(\lk v) =
2h_2(\Delta) + (d-1)h_1(\Delta) \geq 2(d-1)h_1(\Delta)$ by the first part.
Inductively, for $3 \leq j \leq d-1$, we have
\begin{eqnarray*}
j\cdot f_{j-1}(\Delta) &=& \sum_{v \in \Delta}f_{j-2}(\lk v) \\
&\geq& \sum_{v \in \Delta}
\frac{1}{j-1}\left[(2^{j-1}-1)\binom{d-2}{j-2}h_1(\lk v) +
  (d-1)\binom{d-2}{j-2} \right] \\
&\geq& (2^j-2)\binom{d-1}{j-1}h_1(\Delta) + \binom{d-1}{j-1}f_0(\Delta) \\
&=& (2^j-1)\binom{d-1}{j-1}h_1(\Delta) + d\binom{d-1}{j-1}=
\psi_{j-1}(n, d-1).
\end{eqnarray*}
The proof for $j=d$ is similar and is omitted. \endproof

%\begin{eqnarray*}
%d\cdot f_{d-1}(\Delta) &=& \sum_{v \in \Delta}f_{d-2}(\lk v)
%\geq \frac{1}{d-1}\left[(2^{d-1}-2)h_1(\lk v) + 2(d-1) \right] \\
%\geq& (2^{d-1}-2)2h_1(\Delta) + 2f_0(\Delta)
%= (2^d-2)h_1(\Delta) + 2d=\psi_{d-1}(n, d-1),
%\end{eqnarray*}
%and the result follows. \endproof

It is worth remarking that at present we do not know whether the assertion
of Theorem~\ref{balanced} is tight when $n$ is not divisible by $d$. We also
do not know if the stacked cross-polytopal spheres
are the only balanced 2-CM complexes satisfying $2h_2=(d-1)h_1$ when $d$
divides $n$.

In the case when $\Gamma$ is a $2j$-dimensional homology sphere,
the Dehn-Sommerville relations assert that $h_j(\Gamma)=h_{j+1}(\Gamma)$.
If we knew that every balanced 2-CM complex $\Gamma$ of
dimension $2j$ satisfies
$h_j(\Gamma)\leq h_{j+1}(\Gamma)$, a proof similar to that %of the first part
of Theorem~\ref{balanced} would imply that
for a balanced 2-CM complex $\Delta$ of dimension $d-1\geq 2j$,
$$
(j+1) h_{j+1}(\Delta)\geq (d-j) h_j(\Delta).
$$

Finally, we observe that the Lower Bound Theorem \cite{Bar, Kalai}
holds not only for simplicial spheres, but also for triangulations of
connected manifolds, and even normal pseudomanifolds  of dimension at
least two. (The latter result is due to Fogelsanger \cite{Fog}.)
Does Theorem~\ref{balanced} hold for balanced
triangulations of such spaces? Using results of \cite{Nevo-2CM}
and standard tools from rigidity theory, one can show
that any connected pure 3-dimensional simplicial
complex all of whose vertex links are 2-CM, satisfies $h_2\geq h_1$.
The proof analogous to that of Theorem \ref{balanced} then
implies that if $\Delta$ is a balanced triangulation of a
manifold of dimension at least, then $3h_2(\Delta)\geq (d-1)h_1(\Delta)$.
This inequality, however, is weaker than that of Theorem \ref{balanced}.

\bigskip\noindent{\large{\bf Acknowledgements}}
We are grateful to Eran Nevo for introducing us to some of the problems
considered in the paper. Our thanks also go to Christos Athanasiadis,
Eran Nevo, and Ed Swartz for helpful discussions.
Michael Goff's research is partially supported  by a graduate fellowship
from NSF Grant DMS-0801152 and Steven Klee's research
by a graduate fellowship from VIGRE NSF Grant DMS-0354131.

\small{

}

\end{document}